\documentclass{lms}
\usepackage{amsmath,amsfonts,amssymb,url}
\def\classification#1{\def\@class{#1}}
\classification{\null}

\newcommand{\st}{s.t.}

\renewcommand{\theenumi}{\alph{enumi}}
\DeclareFontFamily{OT1}{rsfs}{}
\DeclareFontShape{OT1}{rsfs}{n}{it}{<-> rsfs10}{}
\DeclareMathAlphabet{\mathscr}{OT1}{rsfs}{n}{it}

\DeclareMathOperator{\av}{av}

\DeclareMathOperator{\rad}{rad}

\DeclareMathOperator{\Area}{Area}

\newtheorem{prop}{Proposition}[section]

\newtheorem{lem}[prop]{Lemma}

\newunnumbered{biax}{Bilinear postulate}
\numberwithin{equation}{section}
\title{The parity problem for irreducible cubic forms}
\author{H. A. Helfgott}
\classno{11N32 (primary), 11N35, 11N36 (secondary)}
\begin{document}
\maketitle
\begin{abstract}
Let $f\in \mathbb{Z}\lbrack x,y\rbrack$ be an irreducible homogeneous
polynomial of degree $3$. We show that $f(x,y)$ has an even
number of prime factors as often as an odd number of prime factors.
\end{abstract}
\section{Introduction}
 
Let $f\in \mathbb{Z}\lbrack x,y\rbrack$ be a homogeneous, non-constant
polynomial. Then, it is believed,
\begin{equation}\label{eq:honeysuckle}
\lim_{N\to \infty} \frac{1}{N^2} \sum_{-N\leq x,y\leq N} \mu(f(x,y)) = 0 .
\end{equation}
This conjecture can be traced to Chowla (\cite{Ch}, p. 96). It is closely 
related to the Bunyakovsky--Schinzel conjecture on primes represented by
irreducible polynomials. 

The one-variable analogue of (\ref{eq:honeysuckle}) is classical 
for $\deg f = 1$ and quite hopeless for $\deg f> 1$. We know 
(\ref{eq:honeysuckle}) itself when $\deg f\leq 2$.
(The main ideas of the proof go back to de la Vall\'ee-Poussin
(\cite{DVP1}, \cite{DVP2}); see \cite{Heth}, \S 3.3, for an exposition.)
The problem of proving (\ref{eq:honeysuckle}) when $\deg f\geq 3$ has
remained open until now: 
sieving is forestalled by the parity problem (\cite{Se}), which 
Chowla's conjecture may be said to embody in its pure form.

We prove (\ref{eq:honeysuckle}) for $f$ irreducible of degree $3$. In a 
companion paper \cite{Hered}, we prove (\ref{eq:honeysuckle})
for $f$ reducible of degree $3$.

In \cite{Hered}, we follow Chowla's original formulation,
 using the Liouville function $\lambda$ instead of
$\mu$ in (\ref{eq:honeysuckle}). For $\deg f = 3$, the two formulations
are equivalent: see \S \ref{sec:finrem}.

\subsection{Background}
Problems such as the one considered in this paper were until recently
considered intractable. In 1998, Friedlander and Iwaniec (\cite{FI1}, 
\cite{FI2}) proved that there are infinitely many primes of the form
$x^2 + y^4$. 
One of the main difficulties in adapting this approach to a polynomial 
$f$ other than $x^2 + y^4$ resides in the need to prove that $f$ satisfies
a certain bilinear condition. Another difficulty is that the framework
in \cite{FI2} seemingly breaks down when $f$ represents no more than
$O_\epsilon(N^{2/3+\epsilon})$ integers from $1$ to $N$ for every $\epsilon>0$.

Heath-Brown and Moroz have proved (\cite{HB}, \cite{HBM}, \cite{HBM2})
that any irreducible, homogeneous cubic polynomial $f$ satisfies a bilinear
condition akin to that demanded by \cite{FI2}. Since such an $f$ represents
no more than $O_\epsilon(N^{2/3+\epsilon})$ integers from $1$ to $N$,
Heath-Brown had to abandon the setup in \cite{FI2}, which is based
essentially on Vaughan's identity, for one based on Buchstab's identity.

The framework in \cite{FI2} is quite flexible, and can easily be adapted
to show that $\mu(x^2 + y^4)$ averages to zero.
Unfortunately, it is not clear that \cite{HB} could possibly
 adapted in this way;
Buchstab's identity is in some sense less general than Vaughan's, or to 
a greater extent a statement only about primes.

Our strategy will be to extend the original Friedlander-Iwaniec method
to sequences containing at least $N^{2/3} (\log N)^{-A}$ 
integers from $1$ to $N$, where $A>0$.
This extension goes slightly beyond the natural reach of the method.
Of the technical innovations required, the following may
 be applicable in a wider context.

\subsection{Anti-sieving}
The general situation is as follows. We are given the task of estimating a 
sum $\sum_{a b\leq x} F_{ab}$. Assume we know how to estimate
\begin{equation}\label{eq:afr}
\mathop{\sum_{a b\leq x}}_{a\leq x^{\alpha}/y(x)} F_{ab} \;\;\;\;\;
\text{ and }\;\;\;\;\;
 \mathop{\sum_{a b\leq x}}_{a\geq x^{\alpha} y(x)} F_{ab},\end{equation}
where $0\leq \alpha\leq 1$ and
$\log y(x) = o((\log x)^{1/2})$. 
In order to eliminate the missing interval $x^{\alpha}/y(x)<a<x^{\alpha} y(x)$,
we use a sieve $\lambda_d$ with the primes larger than $y(x)^2$ as our
sieving set:
\[
\mathop{\sum_{a b\leq x}}_{x^{\alpha}/y(x)<a<x^{\alpha} y(x)} F_{a b} =
\mathop{\sum_{a b\leq x}}_{x^{\alpha}/y(x)<a<x^{\alpha} y(x)} 
\left(\sum_{d|a} \lambda_d\right) F_{a b} -
\mathop{\sum_{a b\leq x}}_{x^{\alpha}/y(x)<a<x^{\alpha} y(x)}\;\; 
\mathop{\sum_{d|a}}_{d>(y(x))^2} \lambda_d F_{a b} \]
Notice the use of a sieve as a combinatorial identity, rather than as a
means of approximation; cf. \cite{Hered}, \S 2.1. 
The first term on the right is typically at most
\[O\left(\frac{(\log y(x))^2}{\log x}\right) X,\]
where $X = \sum_{1\leq a\leq x} F_{a}$. The second term on the right equals
\[\mathop{\sum_{a b\leq x}}_{a<x^{\alpha}/y(x)} 
\left(
\mathop{\mathop{
\sum_{d|b}}_{d>\max((y(x))^2,x^{\alpha}/(a y(x)))}}_{d < x^{\alpha}
y(x)/a} \lambda_d\right)
 F_{a b} ,\] which is
akin to the first sum in (\ref{eq:afr}), and can often be treated by
the same methods.

\subsection{Acknowledgements}
Thanks are due to the author's thesis adviser, H. Iwaniec.

\section{Notation}
\subsection{Ideals}
Given a number field $K/\mathbb{Q}$, we write $\mathscr{O}_K$
for the ring of integers of $K$ and $I_K$ for the semigroup of
non-zero ideals of $\mathscr{O}_K$. For $\mathfrak{a}\in I_K$,
define $\tau(\mathfrak{a})$ to be the number of divisors of $\mathfrak{a}$,
$\omega(\mathfrak{a})$ to be the number of prime divisors of $\mathfrak{a}$
and $\mu(\mathfrak{a})$ to be $(-1)^{\omega(\mathfrak{a})}$ if $\mathfrak{a}$
is not divided by the square of any element of $I_K$ (set
 $\mu(\mathfrak{a})=0$ otherwise). Define $\rad(\mathfrak{a}) =
\prod_{\mathfrak{p}|\mathfrak{a}} \mathfrak{p}$.

For $\mathfrak{a}\in I_K$ and $\mathscr{S}$ a finite set of prime ideals
of $\mathscr{O}_K$, we define $r_{\mathscr{S}}(\mathfrak{a}) =
\prod_{\mathfrak{p}\in \mathscr{S}} \mathfrak{p}^{v_{\mathfrak{p}}(\mathfrak{a})}$ and $r_{\backslash \mathscr{S}}(\mathfrak{a}) =
\prod_{\mathfrak{p}\notin \mathscr{S}} \mathfrak{p}^{v_{\mathfrak{p}}(\mathfrak{a})}$.

\subsection{Lattices and convex sets}
A \emph{lattice} is a subgroup of $\mathbb{Z}^2$ of finite index; a \emph{%
lattice coset}
 is a coset of such a subgroup. By the {\em index} $\lbrack \mathbb{Z}^2 : L
\rbrack$ of a lattice
coset $L$ we mean the index of the lattice of which it is a coset.
For any lattice cosets $L_1$, $L_2$ 
with $\gcd(\lbrack \mathbb{Z}^n : L_1 \rbrack , \lbrack \mathbb{Z}^n :
L_2 \rbrack) = 1$, the intersection $L_1\cap L_2$ is a lattice coset with 
\begin{equation}
\lbrack \mathbb{Z}^n : L_1\cap L_2\rbrack = \lbrack \mathbb{Z}^n :
L_1\rbrack \lbrack \mathbb{Z}^n : L_2\rbrack .
\end{equation}
For $L\subset \mathbb{Z}^2$ a lattice coset and $S\subset \mathbb{R}^2$ 
a convex set contained in a square of side $N$,
\begin{equation}  \label{eq:wbrl}
\#(S\cap L) = \frac{\Area(S)}{\lbrack \mathbb{Z}^2 :L\rbrack} + O(N),
\end{equation}
where the implied constant is absolute.

\subsection{Shorthand}
We adopt the following convention from 
\cite{FI2}: given a property $P$, we define
\[f(P(x)) = \begin{cases} f(x) &\text{if $P(x)$ holds,}\\ 0 &\text{otherwise.}
\end{cases}\]
For example, $f(u\leq x\leq v)$ equals $f(x)$ if $u\leq x\leq v$, and $0$ 
otherwise. We will abuse notation by writing $\mathfrak{a}>u$ (or
$\mathfrak{a}<u$, $\mathfrak{a}\leq u$, etc.) when we mean
$N \mathfrak{a}>u$ (or $N \mathfrak{a}<u$, etc.). Thus
 $f(u\leq \mathfrak{a}\leq v)$ equals $f(\mathfrak{a})$ if 
$u\leq N \mathfrak{a}\leq v$, and $0$ otherwise.

We will write $f(x)\ll g(x)$ or $f(x) = O(g(x))$ to mean that 
$|f(x)|\leq c g(x)$ for some positive constant $c$.
By $\mathfrak{p}$ we shall always mean a prime ideal, 
and by $p$ a rational prime.
\section{Postulates}
Let $K/\mathbb{Q}$ be a number field. Let a sequence 
$\{a_{\mathfrak{a}}\}_{\mathfrak{a} \in I_K}$ of non-negative reals 
be given.
We are required to show that
the average of $a_{\mathfrak{a}} \mu(\mathfrak{a})$ over all 
$\mathfrak{a}\in I_K$ is zero.
We need only certain properties of $\{a_{\mathfrak{a}}\}_{\mathfrak{a} \in
I_K}$ to show as much; it will be convenient to list them as postulates for
further reference, and to prove them before the beginning of the main
argument. Weaker postulates might have been used to the detriment of clarity.


As is customary, we define
\[A(n) = \mathop{\sum_{\mathfrak{a}\in I_K}}_{N \mathfrak{a} \leq n} a_{\mathfrak{a}},\;\;\;\;
A_{\mathfrak{d}}(n) = 
\mathop{\sum_{\mathfrak{a}\in I_K}}_{\mathfrak{d}|\mathfrak{a},\; N \mathfrak{a} \leq n} a_{\mathfrak{a}}\]
for any $n>0$, $\mathfrak{d}\in I_K$. 
When bounding
$\sum_{\mathfrak{a} \in I_K : N \mathfrak{a}\leq n} a_{\mathfrak{a}}
\mu(\mathfrak{a})$, we may fix $n$, and require
the postulates only for that particular value of $n$.
\subsection{Statements}
We will work with an approximation
\[A_{\mathfrak{d}}(n) = g(\mathfrak{d}) A(n) + r_{\mathfrak{d}},\]
for some $g:I_K\to \lbrack 0,1)$.
Assume that there is a fixed $D_0\in \mathbb{Z}^+$ and a
$D_1\in \mathbb{Z}^+$ with $\gcd(D_0,D_1)=1$,
$D_1\ll (\log n)^{\rho}$, $\rho>0$, such that the following postulates hold:
\begin{enumerate}
\item\label{it:axprim}
For every prime ideal $\mathfrak{p}\in I_K$ with 
$N \mathfrak{p}$ prime and $\mathfrak{p}\nmid D_0 D_1$, 
\[g(\mathfrak{p}^{\alpha})
 = \frac{(N \mathfrak{p})^{-\alpha}}{1 + 1/ N \mathfrak{p}}
\] for all $\alpha\geq 1$. For $\mathfrak{p}\in I_K$ prime, 
$\mathfrak{p}\nmid D_0$, with $N \mathfrak{p}$ non-prime, we have
$g(\mathfrak{p}^{\alpha}) = 0$ for all $\alpha\geq 1$. 
For $\mathfrak{p} \in I_K$ prime, $\mathfrak{p} | D_1$, with
$N \mathfrak{p}$ prime, either $g(\mathfrak{p}^{\alpha}) = 0$
or $g(\mathfrak{p}^\alpha) = (N \mathfrak{p})^{-\alpha}$
holds for all $\alpha \geq 1$.
\item\label{it:axq1}
For $\mathfrak{a},\mathfrak{b} \in I_K$ with $\gcd(N \mathfrak{a},
N \mathfrak{b}) = 1$,
\[g(\mathfrak{a} \mathfrak{b}) = g(\mathfrak{a}) g(\mathfrak{b}) .\]
\item\label{it:axq2}
For any two distinct 
prime ideals $\mathfrak{p}_1, \mathfrak{p}_2 \nmid D_0$ lying above
the same rational prime,
we have $g(\mathfrak{p}_1 \mathfrak{p}_2 \mathfrak{a}) = 0$ for
all $\mathfrak{a}\in I_K$.
\item\label{it:axr} 
For any $C_1, C_2\geq 0$, there is a $\varkappa>0$ such that
\[\mathop{\sum_{\mathfrak{a}\in I_K}}_{N \mathfrak{a}
\leq n^{2/3} (\log n)^{-\varkappa}} (\tau(\mathfrak{a}))^{C_1}
r_{\mathfrak{a}} \ll_{C_1, C_2} A(n) (\log n)^{-C_2} .\]
\item\label{it:axsq}
For any $C_1, C_2 \geq 0$, there is a $\varkappa>0$ such that
\[
\mathop{\sum_{\mathfrak{d} \in I_K}}_{N \mathfrak{d} > 
(\log n)^{\varkappa}} (\tau(\mathfrak{d}))^{C_1}
 A_{\mathfrak{d}^2}(n) \ll_{C_1,C_2} A(n) (\log n)^{-C_2} .
\]
\item\label{it:axcrude}
For any $C_1, C_2 \geq 0$, there is a $\varkappa>0$ such that
\[\mathop{\sum_{\mathfrak{a}\in I_K}}_{N \mathfrak{a}\leq n (\log n)^{-\varkappa}}
(\tau(\mathfrak{a}))^{C_1} a_{\mathfrak{a}} 
\ll_{C_1,C_2} A(n) (\log n)^{-C_2} .
\]
\item\label{it:biax}
Let $b,c:I_K\to \mathbb{R}$ be any functions with
$b(\mathfrak{a}) \ll (\tau(\mathfrak{a}))^{C_1}$, 
$c(\mathfrak{a}) \ll (\tau(\mathfrak{a}))^{C_2}$.
 Let
\begin{equation}\label{eq:arnos}
d(\mathfrak{a}) = \mathop{\sum_{\mathfrak{d}|\mathfrak{a}}}_{
\gcd(\mathfrak{d},D)=1} c(\mathfrak{a}/{\mathfrak{d}})  
\mu(\mathfrak{d}> \ell) 
,\end{equation}
where $D | D_0 D_1$ and $\ell\gg n^{\epsilon}$
for some $\epsilon>0$. Then, for any $C_3\geq 0$, there is
a $\varkappa>0$ such that
\[\mathop{\mathop{\sum_{\mathfrak{a}, \mathfrak{b}}}_{
\mathfrak{a} \mathfrak{b}\leq n}}_{v \leq N \mathfrak{b} < 2 v}
b(\mathfrak{a}) d(\mathfrak{b}) a_{\mathfrak{a} \mathfrak{b}} 
\ll_{C_1,C_2,C_3} A(n) (\log n)^{-C_3} \]
for all $v\in \lbrack n^{1/2} (\log n)^{\varkappa}, 
n^{3/2} (\log n)^{-\varkappa}\rbrack$. 
Both $\varkappa$ and the implied constant are independent
of $b$, $c$, $D$ and $\ell$.
\end{enumerate}

In brief -- (\ref{it:axprim}) is a statement on $g(\mathfrak{p}^{\alpha})$,
(\ref{it:axq1}) and (\ref{it:axq2}) establish what one may call 
the quasi-multiplicativity of $g$, (\ref{it:axr}) 
states that the residues are small
enough, (\ref{it:axsq}) states that few ideals $\mathfrak{a}$ with
$a_{\mathfrak{a}}\ne 0$ have large square factors,
and (\ref{it:axcrude}) 
is a weak bound on growth. All of these are postulates of a
classical kind (``type I'') whereas (\ref{it:biax}) is a bilinear condition
(and thus of ``type II'').

\subsection{Verification}
Let $K/\mathbb{Q}$ be a cubic extension of $\mathbb{Q}$.
Let $\omega_1, \omega_2\in \mathfrak{O}_K$ be $\mathbb{Q}$--linearly 
independent. 
Let $\varpi>0$ be an arbitrary constant. Given a lattice 
$L\subset \mathbb{Z}^2$, we write $\mathfrak{b}_L$ for the minimal ideal
of $\mathscr{O}_K$ containing the image of $L$ under the map
$(x,y) \mapsto (x \omega_1 + y \omega_2)$.

We must show that, for any sufficiently large $N$, any convex subset
$S \subset \lbrack -N, N\rbrack^2$ with $\Area(S)>N^2 (\log N)^{-\varpi}$,
and any lattice coset $L$ of
index $\lbrack \mathbb{Z}^2 : L\rbrack \leq (\log N)^{\varpi}$, the sequence
\[a_{\mathfrak{a}} = 
\mathop{\mathop{\sum_{(x,y)\in S\cap L}}_{\gcd(x,y)=1}}_{
(x \omega_1 + y \omega_2) = \mathfrak{a} \mathfrak{b}_L} 1\]
satisfies the postulates with $n = \max_{(x,y)\in S} 
|N (x \omega_1 + y \omega_2)| \sim \text{const} \cdot N^3$.
The constants depend only on $K$, $\omega_1$, 
$\omega_2$ and $\varpi$. 
\subsubsection{Linear postulates.}
By \cite{HBM}, Lemma 2.2, there is a $D_0\in \mathbb{Z}^+$
such that, for any prime ideal $\mathfrak{p}\in I_K$, if
$\mathfrak{p}\nmid D_0$ and $\mathfrak{p} | x \omega_1 + y \omega_2$
for some coprime $x,y\in \mathbb{Z}$, then $N \mathfrak{p}$ is prime.
Define $D_1 = \lbrack \mathbb{Z}^2 : L\rbrack$. 

Given $\mathfrak{d}\in I_K$, let $L_{\mathfrak{d}} = \{x,y\in
L : \mathfrak{d} \mathfrak{b}_L | x \omega_1 + y \omega_2\}$. 
For any $\mathfrak{d}\in I_K$ whose norm is a prime power 
$p^{\alpha}$,
define
\[g(\mathfrak{d}) = \frac{\lbrack \mathbb{Z}^2 : L_{\mathfrak{d}}
\rbrack^{-1} - \lbrack \mathbb{Z}^2 : p \mathbb{Z}^2 \cap L_{\mathfrak{d}} 
\rbrack^{-1}}{
\lbrack \mathbb{Z}^2 : L\rbrack^{-1} -
\lbrack \mathbb{Z}^2 : p \mathbb{Z}^2 \cap L\rbrack^{-1}} .\]
(If $L_{\mathfrak{d}} = \emptyset$, we set $g(\mathfrak{d}) = 0$.)
Then postulate (\ref{it:axprim}) follows easily. We may take postulate 
(\ref{it:axq1}) as a definition, and then postulate (\ref{it:axq2}) also follows.
Postulate (\ref{it:axcrude}) is a routinary consequence of (\ref{eq:wbrl})
(cf. \cite{FI1}, p. 1047).

Postulate (\ref{it:axr}) is in essence the same as Lemma 2.2 in
\cite{HBM2}, Lemma 3.2 in \cite{HBM}, or, ultimately, Lemma 5.1 in
\cite{HB}. The contribution of $\mathfrak{a}$
with square factors can be bounded easily (use \cite{Hesq},
Lemma A.5, for factors $\gg (\log N)^C$), and the validity for
all convex sets $S$ can be obtained by partitioning them into
squares\footnote{While \cite{HB} takes as its object a square with
a corner of the form $(x,y)$, $x=y$, its arguments work for any square
whose corners $(x,y)$ satisfy $N (\log N)^{-C} \ll x, y \ll N (\log N)^C$ and whose sides, as in \cite{HB}, are $\gg
N (\log N)^{-C}$.}
of side $\sim N (\log N)^{-C}$. 

The non-trivial part of postulate (\ref{it:axsq}) resides in bounding
the contribution of terms with $\mathfrak{d}$ prime,
$N \mathfrak{d} \geq n (\log n)^{-C}$. Use the bound on the
number of points per fibre in, e.g.,
\cite{Hesq}, Proposition 4.13, together with the crudest bound for
$\tau$. (The techniques in \cite{Grs} would be enough:
postulate (\ref{it:axsq}) is simply stating
that few values of a homogeneous cubic have large square factors.)
\subsubsection{Bilinear postulate.}
While our general framework resembles that of Friedlander and Iwaniec
(\cite{FI2}), the bilinear condition (\ref{it:biax}) is that of
Heath-Brown (\cite{HB}, \cite{HBM}, \cite{HBM2}). 
As (\ref{it:biax}) is an postulate
of the less familiar kind, it is worthwhile to specify the minor changes
we must make to the statement and the proof of Proposition 6.1(ii)
in \cite{HBM}.

First, note that any function defined as in \ref{eq:arnos}
can be extended to ideal numbers so as to fulfil condition (6.1)
in \cite{HBM}, i.e., so as to average to zero at least as fast
as $\exp(-c\sqrt{\log x})$
when restricted to particular ideal classes and lattices of index
$\ll (\log x)^C$. This is simply Siegel-Walfisz; here the condition
$\ell \gg n^{\epsilon}$ in (\ref{it:biax}) is crucial.

While the conditions on $b({\mathfrak{a}})$ are
left unspecified in Proposition 6.1(ii) of \cite{HBM}, it is enough
to have $b(\mathfrak{a}) \ll (\tau(\mathfrak{a}))^{C}$. The term
$b(\mathfrak{a})$ disappears by Cauchy's inequality before the
second equation of p. 279 of \cite{HBM}. While a lacunarity condition
on $b(\mathfrak{a})$ (not fulfilled here) is implicitly used to
eliminate small common factors at the beginning of the proof of
\cite{HBM}, Proposition 6.1 (see also \cite{HB}, \S 11), we may 
remove the condition $\gcd(x,y)=1$ by carrying the argument in 
\cite{HBM}, pp. 278--284, for each
lattice coset $L_d$ defined by
$(x,y)\in L, d|\gcd(x,y)$, where $d\ll (\log N)^C$, and then sieving
out such lattice cosets.
A riddle (vd. Proposition 3.2 and Corollary 3.3 in \cite{Hesq}) 
enables
us to sieve with a total error 
term of relative size $O((\log N)^{-C + O(1)})$
while handling only the $L_d$'s with $d\ll (\log N)^C$;
cf. \cite{Hesq}, Propositions 3.11 and 3.12.
\section{Proof}\label{sec:ufeuf}
\begin{lem}[(A variant of Vaughan's identity)]\label{lem:ontont}
Let $y$, $u$, $w$ be positive numbers satisfying
 $y\leq u\leq w$. Let $K/\mathbb{Q}$ be a number field; let
$\mathscr{Q}$ be a finite set of prime ideals thereof, and 
let $h:I_K\to \mathbb{C}$ be an arbitrary function.
Then, for any $\mathfrak{a}\in I_K$,
\[h(\mathfrak{a})
 = \sum_{j=1}^4 \beta_j(\mathfrak{a}) -
 \sum_{j=5}^7 \beta_j(\mathfrak{a}) ,\]
where
\begin{equation}\label{eq:sensat}\begin{aligned}
\beta_1(\mathfrak{a}) &= h(\mathfrak{a}\leq u) +
 \sum_{*} h(\mathfrak{b}) \mu(\mathfrak{c}
\leq u),\\
\beta_2(\mathfrak{a}) &= 
\sum_* h(u<\mathfrak{b}\leq w) \mu(\mathfrak{c}>u),\;\;\;
\beta_3(\mathfrak{a}) = \sum_* h(\mathfrak{b}> w) 
\mu(u<\mathfrak{c}\leq w),\\
\beta_4(\mathfrak{a}) &= \sum_*
h(\mathfrak{b}> w) \mu(\mathfrak{c}>w),\;\;\;
\;\;\;\;\;
\;\;\beta_5(\mathfrak{a}) = \sum_* 
h(\mathfrak{b} \leq u) \mu(\mathfrak{c} \leq y),\\
\beta_6(\mathfrak{a}) &= \sum_* 
h(\mathfrak{b} \leq y) \mu(y<\mathfrak{c} \leq u),\;\;\;\;
\beta_7(\mathfrak{a}) = \
\sum_* 
h(y<\mathfrak{b} \leq u) \mu(y<\mathfrak{c} 
\leq u),
\end{aligned}\end{equation}
and $\sum_*$ stands for 
$\sum_{\mathfrak{b} \mathfrak{c} | \mathfrak{a},\;
r_{\mathscr{Q}}(\mathfrak{b}) = r_{\mathscr{Q}}(\mathfrak{a})}$.
\end{lem}
\begin{proof}
For any $\mathfrak{a}\in I_K$, 
\[
\beta_2(\mathfrak{a}) + \beta_3(\mathfrak{a}) + \beta_4(\mathfrak{a}) =
\sum_* h(\mathfrak{b}>u) \mu(\mathfrak{c}>u),\;\;\;
\beta_5(\mathfrak{a}) + \beta_6(\mathfrak{a}) + \beta_7(\mathfrak{a}) =
\sum_* h(\mathfrak{b}\leq u) \mu(\mathfrak{c}
\leq u) .\]
Thus
\[\begin{aligned}
h(\mathfrak{a}) &= h(\mathfrak{a}\leq u) + h(\mathfrak{a}>u) =
h(\mathfrak{a}\leq u) + \sum_* h(\mathfrak{b}>u) \mu(
\mathfrak{c}) 
= h(\mathfrak{a}\leq u) + \sum_* h(\mathfrak{b}) \mu(
\mathfrak{c}\leq u) \\ &- \sum_* h(\mathfrak{b}\leq u) 
\mu(\mathfrak{c}\leq u) 
+ \sum_* h(\mathfrak{b}>u) \mu(\mathfrak{c}>u) 
= \beta_1(\mathfrak{a}) - \sum_{j=5}^7 \beta_j(\mathfrak{a}) +
\sum_{j=2}^4 \beta_j(\mathfrak{a}) .\end{aligned}\]
\end{proof}

Set $h=\mu$ in Lemma \ref{lem:ontont}.
Let $z = e^{(\log \log x) (\log \log \log x)^{\epsilon/2}}$ for
some $\epsilon>0$. Let
$y = x^{1/3} z^{-2}$, $u = x^{1/3} z$, $w = x^{1/2} z^{-1}$,
$\mathscr{Q} = \{\mathfrak{p} : \mathfrak{p}|D_0 D_1\}$.
Since we wish to estimate 
$\sum_{\mathfrak{a}\leq x} 
a_{\mathfrak{a}} \mu(\mathfrak{a})$, 
we will evaluate $\sum_{\mathfrak{a}\leq x} 
a_{\mathfrak{a}} \beta_j(\mathfrak{a})$ for $j=1,2,\dotsc,7$.
The cases $j=1$, $j=5$ and $j=6$ are easy.
\begin{lem}\label{lem:j1}
Let $\{a_{\mathfrak{a}}\}$ be a sequence satisfying the postulates for
$n=x$. 
Then,
for any $C>0$,
\[\sum_{\mathfrak{a}\leq x} a_{\mathfrak{a}} \beta_1(\mathfrak{a})
\ll A(x) (\log x)^{-C}.\] 
\end{lem}
\begin{proof}
By postulate (\ref{it:axcrude}), we have $\sum_{\mathfrak{a}\leq x} 
a_{\mathfrak{a}} \mu(\mathfrak{a}\leq u)\ll A(x) (\log x)^{-C}$ and
\[\sum_{\mathfrak{a}\leq x} \sum_* a_{\mathfrak{a}} 
\mu(\mathfrak{b}) \mu(\mathfrak{c}\leq u)  =
\sum_{\mathfrak{a}\leq x} a_{\mathfrak{a}} \mu(r_{\mathscr{Q}}(\mathfrak{a}))
\mu(r_{\backslash \mathscr{Q}}(\mathfrak{a}) \leq u) \,\ll
\sum_{\mathfrak{a}\leq u D_0^3 D_1^3} \tau(\mathfrak{a}) 
a_{\mathfrak{a}} \ll A(x) (\log x)^{-C}.\]
\end{proof}
\begin{lem}\label{lem:j5}
Let $\{a_{\mathfrak{a}}\}$ be a sequence satisfying the postulates for
$n=x$. 
 Then,
for any $C>0$,
\[\sum_{\mathfrak{a}\leq x} a_{\mathfrak{a}} \beta_5(\mathfrak{a})
\ll_C A(x) (\log x)^{-C} .\] 
\end{lem}
\begin{proof}
By exclusion--inclusion and postulate (\ref{it:axr}), 
\[\begin{aligned}
\sum_{\mathfrak{a}\leq x} a_{\mathfrak{a}} \beta_5(\mathfrak{a})
&= \sum_{\mathfrak{e} : r_{\mathscr{Q}}(\mathfrak{e}) = \mathfrak{e}}
\mu(\mathfrak{e})
\mathop{\sum_{\mathfrak{b}\leq u}}_{r_{\backslash \mathscr{Q}}(\mathfrak{b})
= \mathfrak{b}} \mu(\mathfrak{b}) 
\mathop{\sum_{\mathfrak{c}\leq y}}_{r_{\backslash \mathscr{Q}}(\mathfrak{c})
= \mathfrak{c}} \mu(\mathfrak{c})
\mathop{\sum_{\mathfrak{d}\leq x/ N (\mathfrak{e} \mathfrak{b} \mathfrak{c})}
}_{r_{\backslash \mathscr{Q}}(\mathfrak{d}) =
\mathfrak{d}}
a_{\mathfrak{e} \mathfrak{b} \mathfrak{c} \mathfrak{d}}\\
&= \sum_{\mathfrak{e} : r_{\mathscr{Q}}(\mathfrak{e}) = \mathfrak{e}}
\mu(\mathfrak{e})
\sum_{\mathfrak{f} : r_{\mathscr{Q}}(\mathfrak{f}) = \mathfrak{f}}
\mu(\mathfrak{f})
\mathop{\sum_{\mathfrak{b}\leq u}}_{r_{\backslash \mathscr{Q}}(\mathfrak{b})
= \mathfrak{b}} \mu(\mathfrak{b}) 
\mathop{\sum_{\mathfrak{c}\leq y}}_{r_{\backslash \mathscr{Q}}(\mathfrak{c})
= \mathfrak{c}} \mu(\mathfrak{c})
\sum_{\mathfrak{d}\leq x/ N(\mathfrak{e} \mathfrak{f} \mathfrak{b}
\mathfrak{c})}
a_{\mathfrak{e} \mathfrak{f} \mathfrak{b} \mathfrak{c} \mathfrak{d}} \\
&\leq 
 A(x)
\sum_{\mathfrak{e} : r_{\mathscr{Q}}(\mathfrak{e}) = \mathfrak{e}}
\mu(\mathfrak{e})
\sum_{\mathfrak{f} : r_{\mathscr{Q}}(\mathfrak{f}) = \mathfrak{f}}
\mu(\mathfrak{f})
\mathop{\sum_{\mathfrak{b}\leq u}}_{r_{\backslash \mathscr{Q}}(\mathfrak{b})
= \mathfrak{b}} \mu(\mathfrak{b}) g(\mathfrak{b})
\mathop{
\mathop{\sum_{\mathfrak{c}\leq y}}_{
r_{\backslash \mathscr{Q}}(\mathfrak{c})
= \mathfrak{c}}}_{\gcd(N \mathfrak{b}, N \mathfrak{c}) = 1}
 \mu(\mathfrak{c}) g(\mathfrak{c}) \\
&+ O(A(x) (\log x)^{-C}) .\end{aligned}\]
By the standard zero-free region for $\zeta_K(s)$, the innermost sum
$\sum_{\mathfrak{c}}$ is $\ll  e^{-c \sqrt{\log x}}$, $c>0$.
\end{proof}
\begin{lem}\label{lem:j6}
Let $\{a_{\mathfrak{a}}\}$ be a sequence satisfying the postulates for
$n=x$. 
Then,
for any $C>0$,
\[\sum_{\mathfrak{a}\leq x} a_{\mathfrak{a}} \beta_6(\mathfrak{a})
\ll_C A(x) (\log x)^{-C}.\] 
\end{lem}
\begin{proof}
Same as Lemma \ref{lem:j5}.
\end{proof}
The following two lemmas are direct consequences of the bilinear postulate.
\begin{lem}\label{lem:j2}
Let $\{a_{\mathfrak{a}}\}$ be a sequence satisfying the postulates for
$n=x$. 
Then, for any $C>0$,
\[\sum_{\mathfrak{a}\leq x} a_{\mathfrak{a}} \beta_2(\mathfrak{a})
\ll_C A(x) (\log x)^{-C} .\] 
\end{lem}
\begin{proof}
Apply postulate (\ref{it:biax}) with $D = D_0 D_1$, $\ell = u$, 
\[b(\mathfrak{a}) = \mu(u<\mathfrak{a}\leq w),\;\;\;\;\;
c(\mathfrak{a}) = \begin{cases}
1 &\text{if $\gcd(\mathfrak{a},D_0 D_1) = 1$,}\\
0 &\text{otherwise,}\end{cases}\]
and $v$ ranging from $x w^{-1} (\log x)^{-C}$ to $x u^{-1}/2$. Use
postulate (\ref{it:biax}) to bound the remaining terms.
\end{proof}
\begin{lem}\label{lem:j3}
Let $\{a_{\mathfrak{a}}\}$ be a sequence satisfying the postulates for
$n=x$. 
Then,
for any $C>0$,
\[\sum_{\mathfrak{a}\leq x} a_{\mathfrak{a}} \beta_3(\mathfrak{a})
\ll_C A(x) (\log x)^{-C} .\] 
\end{lem}
\begin{proof}
Apply postulate (\ref{it:biax}) with $D = 1$, $\ell = w$,
\[b(\mathfrak{a}) = \begin{cases}
\mu(u<\mathfrak{a}\leq w) &\text{if $\gcd(\mathfrak{a},D_0 D_1) = 1$,}\\
0 &\text{otherwise,}\end{cases}
\;\;\;\;\;
c(\mathfrak{a}) = \begin{cases}
1 &\text{if $\gcd(\mathfrak{a},D_0 D_1) = 1$,}\\
0 &\text{otherwise,}\end{cases}\]
and $v$ ranging from $x w^{-1} (\log x)^{-C}$ to $x u^{-1}/2$. Use
postulate (\ref{it:biax}) to bound the remaining terms.
\end{proof}
It remains to consider the sums
$\sum_{\mathfrak{a}} a_{\mathfrak{a}} \beta_j(\mathfrak{a})$
for $j=4,7$. We will recur to anti-sieving and a certain kind
of cancellation.
\begin{lem}\label{lem:j7}
Let $\{a_{\mathfrak{a}}\}$ be a sequence satisfying the postulates for
$n=x$. 
Then
\[\sum_{\mathfrak{a}\leq x} a_{\mathfrak{a}} \beta_7(\mathfrak{a})
\ll \frac{(\log \log x)^4 (\log \log \log x)^{\epsilon}}{\log x} A(x) ,\]
where $\epsilon$ is as in the definition of $z$.
\end{lem}
\begin{proof}
Let $\{\lambda_{\mathfrak{d}}\}$ be a Rosser-Iwaniec sieve 
\cite{Col}
with sieved set $\mathscr{A}=\{\mathfrak{b}\in I_K:
y<N \mathfrak{b} \leq u\}$, 
multiplicities $\theta(\mathfrak{b}) =  g(\mathfrak{b}) \mu^2(\mathfrak{b})$,
sieving set $\mathscr{P} = \{\mathfrak{p}: u y^{-1} < N \mathfrak{p} 
\leq
w u^{-1}\}$ and upper cut $w u^{-1}$. (Brun's pure sieve would do almost
as well.)
By definition,
$\lambda_{\mathfrak{d}}=0$ if 
$1< N \mathfrak{d} \leq u y^{-1}$ or $N \mathfrak{d} > w u^{-1}$. Since
$\lambda_{1}=1$, we have the identity
\[1 = \sum_{\mathfrak{d}|\mathfrak{b}} \lambda_{\mathfrak{d}} \;\; -
\mathop{\sum_{u y^{-1} < \mathfrak{d} \leq w u^{-1}}}_{\mathfrak{d} | 
\mathfrak{b}} \lambda_{\mathfrak{d}}\]
for every $\mathfrak{b}\in I_K$. Hence
$\beta_7(\mathfrak{a}) = \beta_8(\mathfrak{a}) - \beta_9(\mathfrak{a})$,
where
\[\begin{aligned}
\beta_8(\mathfrak{a})
&= \mathop{\sum_{\mathfrak{b} \mathfrak{c} | \mathfrak{a}}}_{
\mathfrak{p} | \mathfrak{a}/\mathfrak{b} \Rightarrow \mathfrak{p}\notin
\mathscr{Q}} \sum_{\mathfrak{d}|\mathfrak{b}} \lambda_{\mathfrak{d}}
\mu(y<\mathfrak{b} \leq u) \mu(y<\mathfrak{c} \leq u) \\
\beta_9(\mathfrak{a}) &=
\mathop{\sum_{\mathfrak{b} \mathfrak{c} | \mathfrak{a}}}_{
\mathfrak{p} | \mathfrak{a}/\mathfrak{b} \Rightarrow \mathfrak{p}\notin
\mathscr{Q}} \mathop{\sum_{u y^{-1} < \mathfrak{d} \leq w u^{-1}}}_{
\mathfrak{d}|\mathfrak{b}} \lambda_{\mathfrak{d}}
\mu(y<\mathfrak{b} \leq u) \mu(y<\mathfrak{c} \leq u) .
\end{aligned}\]
We begin by bounding $\sum_{\mathfrak{a}\leq x} a_{\mathfrak{a}}
\beta_9(\mathfrak{a})$. Changing the order of summation,
we obtain
\[\begin{aligned}
\sum_{\mathfrak{a}\leq x} a_{\mathfrak{a}} \beta_9(\mathfrak{a})
&= \mathop{\sum_{y < \mathfrak{c} \leq u}}_{\mathfrak{p}|
\mathfrak{c} \Rightarrow \mathfrak{p}\notin \mathscr{Q}}
 \mu(\mathfrak{c})
\sum_{u y^{-1} < \mathfrak{d}\leq w u^{-1}}
\lambda_{\mathfrak{d}}
\sum_{y/N \mathfrak{d} < \mathfrak{e} \leq u/N \mathfrak{d}}
\mu(\mathfrak{d} \mathfrak{e})
\mathop{\sum_{\mathfrak{f} \leq x/ N(\mathfrak{c} \mathfrak{d}
\mathfrak{e})}}_{\mathfrak{p}|\mathfrak{f} \Rightarrow \mathfrak{p}
\notin \mathscr{Q}} a_{\mathfrak{c} \mathfrak{d}
\mathfrak{e} \mathfrak{f}}\\
&= \sum_{u<\mathfrak{g}\leq w} 
\mathop{\mathop{\mathop{\sum_{\mathfrak{d}|\mathfrak{g}}}_{
\mathfrak{d} > u y^{-1},\; \mathfrak{d}\geq N\mathfrak{g}/u}}_{
\mathfrak{d} \leq w u^{-1},\; \mathfrak{d} < N\mathfrak{g}/y}}_{
\mathfrak{p}|\mathfrak{g}/\mathfrak{d}\Rightarrow \mathfrak{p}\notin \mathscr{Q}}
\lambda_{\mathfrak{d}} \mu(\mathfrak{d}) \mu(\mathfrak{g}/\mathfrak{d})
\mathop{\sum_{y/N\mathfrak{d} < \mathfrak{e} \leq u/N 
\mathfrak{d}}}_{\gcd(\mathfrak{e},\mathfrak{d})=1}
 \mu(\mathfrak{e})
\mathop{\sum_{\mathfrak{f} \leq x/ N(\mathfrak{g}
\mathfrak{e})}}_{\mathfrak{p}|\mathfrak{f} \Rightarrow \mathfrak{p}
\notin \mathscr{Q}} a_{\mathfrak{g} \mathfrak{e} \mathfrak{f}}.\end{aligned}\]

Since $\mathfrak{d}$ has no factors of norm less than $z^3$ when
$\lambda_{\mathfrak{d}}\ne 0$, we may remove
the condition $\gcd(\mathfrak{e},\mathfrak{d})$ with an error
of at most $O(A(x) (\log x)^{-C})$ by means of postulates 
(\ref{it:axprim})--(\ref{it:axr}).
We can make the intervals
of summation of $\mathfrak{d}$ and $\mathfrak{e}$ independent
from each other by slicing $\lbrack u y^{-1}, w u^{-1}\rbrack$
into intervals of the form $\lbrack K, K (1 + (\log x)^{-C}))$.
There are at most $O((\log x)^{C+1})$ such intervals, 
and the error incurred
during the slicing is at most $O(A(x) (\log x)^{-C + O(1)})$.
Hence
\[\sum_{\mathfrak{a}\leq x} a_{\mathfrak{a}} \beta_9(\mathfrak{a})
\ll A(x) (\log x)^{-C+O(1)} + (\log x)^{C+1}
\max_{u y^{-1} < K \leq w u^{-1}} \left|
\mathop{\mathop{\sum_{\mathfrak{g}, \mathfrak{h}}}_{
\mathfrak{g} \mathfrak{h} \leq x}}_{u<\mathfrak{g}\leq w} 
b_K(\mathfrak{g}) d_K(\mathfrak{h}) a_{\mathfrak{g} \mathfrak{h}}
\right| ,\]
where
\begin{equation}\label{eq:poeque}
b_K(\mathfrak{g}) = 
\mathop{\mathop{\mathop{\sum_{\mathfrak{d}|\mathfrak{g}}}_{
K\leq \mathfrak{d} < (1 + (\log x)^{-C}) K}}_{
y< \mathfrak{g}/\mathfrak{d}\leq u}}_{
\mathfrak{p}|\mathfrak{g}/
\mathfrak{d} \Rightarrow \mathfrak{p} \notin
\mathscr{Q}} \lambda_{\mathfrak{d}} \mu(\mathfrak{d})
\mu(\mathfrak{g}/\mathfrak{d}),\;\;\;\;\;\;\;\;\;\;
d_K(\mathfrak{h}) = 
\mathop{\mathop{\sum_{\mathfrak{e}|\mathfrak{h}}}_{
y/K<\mathfrak{e}\leq u/K}}_{\mathfrak{p}|\mathfrak{h}/
\mathfrak{e} \Rightarrow \mathfrak{p}\notin \mathscr{Q}}
 \mu(\mathfrak{e}) .\end{equation}
We apply postulate (\ref{it:biax}) with $D=1$, $b = b_K$,
$c(\mathfrak{a})=1$ when $\gcd(\mathfrak{a}, D_0 D_1) = 1$,
$c(\mathfrak{a})=0$ otherwise, and $\ell = y/K$ or $\ell = u/K$
(in succession).
We obtain a bound $\sum_{\mathfrak{a}\leq x} a_{\mathfrak{a}}
\beta_9(\mathfrak{a}) \ll A(x) (\log x)^{-C'}$ for $C'$ arbitrarily
large.

We must now bound 
\begin{equation}\label{eq:hab}
\sum_{\mathfrak{a}\leq x} a_{\mathfrak{a}}
\beta_8(\mathfrak{a})
= 
\sum_{y<\mathfrak{b}\leq u}
\sum_{\mathfrak{d}|\mathfrak{b}}
\lambda_{\mathfrak{d}} \mu(\mathfrak{b})
\mathop{\sum_{\mathfrak{e}\leq x/ N \mathfrak{b}}}_{
\mathfrak{p} |\mathfrak{e} \Rightarrow \mathfrak{p} \notin
\mathscr{Q}}
\sum_{\mathfrak{c}|\mathfrak{e}}
\mu(y<\mathfrak{c}\leq u)
a_{\mathfrak{b} \mathfrak{e}} .\end{equation}
We must find cancellation in the innermost sum
and lower $\mathfrak{b} \mathfrak{c}$ below
$(\log x)^{-\varkappa} x^{2/3}$. For $\mathfrak{e}>1$,
\begin{equation}\label{eq:rust}\begin{aligned}
\sum_{\mathfrak{c}|\mathfrak{e}} \mu(y<\mathfrak{c}\leq u)
&= \sum_{\mathfrak{c}|\mathfrak{e}} \mu(\mathfrak{c})
-
\sum_{\mathfrak{c}|\mathfrak{e}} \mu(\mathfrak{c}\leq y) -
\sum_{\mathfrak{c}|\mathfrak{e}} \mu(\mathfrak{c}>u) \\
&=
-
\sum_{\mathfrak{c}|\mathfrak{e}} \mu(\mathfrak{c}\leq y) 
-
\mu(\rad(\mathfrak{e})) 
\sum_{\mathfrak{c}|\mathfrak{e}} 
      \mu(\mathfrak{c}<\rad(\mathfrak{e})/u) .
\end{aligned}\end{equation}
We will first bound the contribution from
$\sum_{\mathfrak{c}|\mathfrak{e}} \mu(\mathfrak{c}\leq y)$.
A bound for $
\mu(\rad(\mathfrak{e})) 
\sum_{\mathfrak{c}|\mathfrak{e}} 
      \mu(\mathfrak{c}<\rad(\mathfrak{e})/u)$ 
will be obtained later in a similar fashion. The terms with
 $\mathfrak{e}=1$ may be ignored by postulate (\ref{it:axcrude}).

Suppose $\mathfrak{e}$ has a prime divisor 
$\mathfrak{p}\leq l$, where $l>0$ is fixed. Then the set
of all square-free divisors of $\mathfrak{e}$ can be partitioned
into pairs $\{\mathfrak{o},\mathfrak{o} \mathfrak{p}\}$. 
Evidently,
$\mu(\mathfrak{o})= -\mu(\mathfrak{o} \mathfrak{p})$.
We have either $\mathfrak{o} \leq y$,
$\mathfrak{o} \mathfrak{p}\leq y$ or
$\mathfrak{o} > y$, $\mathfrak{o} \mathfrak{p} >y$, unless
$\mathfrak{o}$ lies in the range
$y/l < \mathfrak{c} \leq y$. Hence
\begin{equation}\label{eq:aelita}
\left| \sum_{\mathfrak{c}|\mathfrak{e}} \mu(\mathfrak{c} \leq y)
\right| \leq \mathop{\sum_{\mathfrak{c}|\mathfrak{e}}}_{
y/l<\mathfrak{c}\leq y} 1 .\end{equation}
Define $l_j = 2^{2^j}$ for $j\geq 0$. Note that
$x^{1/2} < l_{\lfloor \log_2 \log_2 x\rfloor} \leq x$.
Let \[\begin{aligned}
L_0 &= \{\mathfrak{e}: 2|\mathfrak{e}, 
(\mathfrak{p}|\mathfrak{e} \Rightarrow \mathfrak{p}\notin
\mathscr{Q})\},\\
L_j &= \{\mathfrak{e}\in I_K :
(\text{$\exists \mathfrak{p}\leq l_j$ \st\; $\mathfrak{p}|\mathfrak{e}$})
\wedge (\forall \mathfrak{p}\leq l_{j-1},\; \mathfrak{p}
\nmid \mathfrak{e}) \wedge
(\mathfrak{p}|\mathfrak{e} \Rightarrow \mathfrak{p}\notin
\mathscr{Q})\} .\end{aligned}\]
Then, by (\ref{eq:aelita}),
\begin{equation}\label{eq:nicom}
\mathop{\sum_{\mathfrak{e}\leq x/N \mathfrak{b}}}_{
\mathfrak{e}\in L_j} \left|\sum_{\mathfrak{c}|\mathfrak{e}}
\mu(\mathfrak{c} \leq y)\right| a_{\mathfrak{b} \mathfrak{e}}
\leq
\mathop{\sum_{\mathfrak{e}\leq x/N \mathfrak{b}}}_{
\mathfrak{e}\in L_j} 
\mathop{\sum_{\mathfrak{c}|\mathfrak{e}}}_{
y/l_j< \mathfrak{c}\leq y} 
a_{\mathfrak{b} \mathfrak{e}} \leq
\mathop{\sum_{y/l_j<\mathfrak{c}\leq y}}_{\mathfrak{p}|
\mathfrak{c}\Rightarrow
\mathfrak{p}>l_{j-1}, \mathfrak{p}\notin \mathscr{Q}}
\sum_{\mathfrak{j}\leq x/N(\mathfrak{b} \mathfrak{c})}
a_{\mathfrak{b} \mathfrak{c} \mathfrak{j}} .\end{equation}
By postulates (\ref{it:axprim})--(\ref{it:axr}), it follows that the
contribution of the sums $\sum_{\mathfrak{c}|\mathfrak{e}}
\mu(\mathfrak{c}\leq y)$ to (\ref{eq:hab}) is at most
\begin{equation}\label{eq:blod}
A(x)
\mathop{\sum_{y<\mathfrak{b}\leq u}}_{\mu(\mathfrak{b}) = \pm 1}
\sum_{\mathfrak{d}|\mathfrak{b}}
\lambda_{\mathfrak{d}} \sum_{j=0}^{\lfloor \log_2 \log_2 x\rfloor}
\mathop{\sum_{y/l_j<\mathfrak{c}\leq y}}_{\mathfrak{p}|
\mathfrak{c}\Rightarrow
\mathfrak{p}>l_{j-1}, \mathfrak{p}\notin \mathscr{Q}}
g(\mathfrak{b} \mathfrak{c}) + O(A(x) (\log x)^{-C}) \end{equation}
in absolute value.
Now $g(\mathfrak{b} \mathfrak{c}) \leq g(\mathfrak{b}) / N \mathfrak{c}$,
and
\[
\mathop{\sum_{y/l_j<\mathfrak{c}\leq y}}_{\mathfrak{p}|
\mathfrak{c}\Rightarrow
\mathfrak{p}>l_{j-1}, \mathfrak{p}\notin \mathscr{Q}}
\frac{1}{N \mathfrak{c}} \ll
  (\log y - \log \frac{y}{l_j}) \cdot
\prod_{\mathfrak{p}\leq l_{j-1}}
 \left(1 - \frac{1}{N \mathfrak{p}}
\right) \ll \frac{\log l_j}{\log l_{j-1}} \ll 1 .\]
Thus, the quadruple sum in (\ref{eq:blod}) is at most
$O(\sum_{y<\mathfrak{b}\leq u, \mu(\mathfrak{b}) = 
\pm 1} \sum_{\mathfrak{d}|\mathfrak{b}}
\lambda_{\mathfrak{d}} g(\mathfrak{b}))\cdot \log \log x$.
The main result on the Rosser-Iwaniec sieve (\cite{Col},
Lemma 3), granted postulates (\ref{it:axprim})--(\ref{it:axq2}), 
gives:
\begin{equation}\label{eq:maleta}\begin{aligned}
\mathop{\sum_{y<\mathfrak{b}\leq u}}_{\mu(\mathfrak{b}) = \pm 1}
 \sum_{\mathfrak{d}|\mathfrak{b}}
\lambda_{\mathfrak{d}} g(\mathfrak{b}) &\ll
\mathop{
\sum_{y<\mathfrak{b}\leq u}}_{
u y^{-1} < \mathfrak{p} < w u^{-1} \Rightarrow \mathfrak{p}\nmid \mathfrak{b}}
|\mu(\mathfrak{b})|
g(\mathfrak{b}) \ll (\log u - \log y) \frac{\log u y^{-1}}{\log w u^{-1}}\\
&\ll \frac{(\log z)^2}{\log x} = \frac{(\log \log x)^2 (\log \log \log x)^{
\epsilon}}{\log x} .\end{aligned}\end{equation}
It remains to bound the contribution from the terms
$\mu(\rad(\mathfrak{e})) 
\sum_{\mathfrak{c}|\mathfrak{e}} 
      \mu(\mathfrak{c}<\rad(\mathfrak{e})/u)$, viz.,
\begin{equation}\label{eq:tugtug}
\mathop{\sum_{y<\mathfrak{b}\leq u}}_{\mu(\mathfrak{b}) = \pm 1}
\sum_{\mathfrak{d}|\mathfrak{b}} \lambda_{\mathfrak{d}} 
\mathop{\sum_{\mathfrak{e}\leq x/N \mathfrak{b}}}_{
\mathfrak{p}|\mathfrak{e} \Rightarrow  \mathfrak{p} \notin \mathscr{Q}}
\left|
\sum_{\mathfrak{c}|\mathfrak{e}} \mu(\mathfrak{c} < \rad(\mathfrak{e})/u) 
\right| a_{\mathfrak{b} \mathfrak{e}}.\end{equation}
By postulates (\ref{it:axsq}) and (\ref{it:axcrude}), there is a $C>0$ such
that the terms with $N \mathfrak{e} \leq x / (\log x)^C N \mathfrak{b}$
contribute less than $A(x)/\log x$ to the total. Proceeding as in
(\ref{eq:aelita})--(\ref{eq:blod}), we obtain that (\ref{eq:tugtug}) is
at most
\[A(x) \mathop{\sum_{y<\mathfrak{b}\leq u}}_{\mu(\mathfrak{b}) = \pm 1}
\sum_{\mathfrak{d} | \mathfrak{b}} \lambda_{\mathfrak{d}}
\sum_{j=0}^{\lfloor \log_2 \log_2 x\rfloor}
\sum_{k=0}^{\lfloor \log_2 \log_2 x\rfloor} \mathop{\sum_{
x/(l_j u 2^{k+1} \cdot N \mathfrak{b}) \leq \mathfrak{c} <
x/(u 2^k \cdot N \mathfrak{b})}}_{\mathfrak{p} | \mathfrak{c} \Rightarrow
\mathfrak{p} > l_{j-1}, \mathfrak{p}\in \mathscr{Q}} g(\mathfrak{b}
\mathfrak{c}) + O(A(x) (\log x)^{-C}) .\]
Again as before, we obtain that this is at most
$O(\sum_{y<\mathfrak{b} \leq u, \mu(\mathfrak{b}) = \pm 1}
\sum_{\mathfrak{d} | \mathfrak{b}} \lambda_{\mathfrak{d}} g(\mathfrak{b}))
A(x) (\log \log x)^2$. By (\ref{eq:maleta}), we are done.
\end{proof}
\begin{lem}\label{lem:j4}
 Let $\{a_{\mathfrak{a}}\}$ be a sequence satisfying the postulates for
$n=x$. Then
\[\sum_{\mathfrak{a}\leq x} a_{\mathfrak{a}} \beta_4(\mathfrak{a})
\ll \frac{(\log \log x)^4 (\log \log \log x)^{\epsilon}}{\log x} A(x) ,\]
where $\epsilon$ is as in the definition of $z$.
\end{lem}
\begin{proof}
Let $\{\lambda_{\mathfrak{d}}\}$ be a Rosser-Iwaniec sieve with sieved
set $\mathscr{A} = \{\mathfrak{b} \in I_K : w < N \mathfrak{b} <
x/w\}$, multiplicities $\theta(\mathfrak{b}) = g(\mathfrak{b})
\mu^2(\mathfrak{b})$, sieving set $\mathscr{P} = \{\mathfrak{p} :
x w^{-2} < N \mathfrak{p} \leq w x^{-1} y^2\}$ and upper cut
$w x^{-1} y^2$. By definition, $\lambda_{\mathfrak{d}} = 0$,
if $1 < N \mathfrak{d} \leq x w^{-2}$ or $N \mathfrak{d} > w x^{-1} y^2$. 
Then
$\beta_{4}(\mathfrak{a}) = \beta_{10}(\mathfrak{a}) - \beta_{11}(\mathfrak{a})$
for $\mathfrak{a}\in I_K$ with $\mathfrak{a}\leq x$, where
\[\begin{aligned}
\beta_{10}(\mathfrak{a}) &= \mathop{\sum_{\mathfrak{b} \mathfrak{c} |
\mathfrak{a}}}_{\mathfrak{p} | \mathfrak{a}/\mathfrak{b} \Rightarrow
\mathfrak{p}\notin \mathscr{Q}} \sum_{\mathfrak{d} | \mathfrak{b}}
\lambda_{\mathfrak{d}} \mu(w < \mathfrak{b} \leq x w^{-1})
\mu(w < \mathfrak{c} \leq x w^{-1})\\
\beta_{11}(\mathfrak{a}) &= \mathop{\sum_{\mathfrak{b} \mathfrak{c} |
\mathfrak{a}}}_{\mathfrak{p} | \mathfrak{a}/\mathfrak{b} \Rightarrow
\mathfrak{p}\notin \mathscr{Q}} 
\mathop{\sum_{x w^{-2} < \mathfrak{d} < w x^{-1} y^2}}_{
\mathfrak{d} | \mathfrak{b}}
\lambda_{\mathfrak{d}} \mu(w < \mathfrak{b} \leq x w^{-1})
\mu(w < \mathfrak{c} \leq x w^{-1}) .\end{aligned}\]
Much as for $\sum a_{\mathfrak{a}} \beta_9(\mathfrak{a})$ in
the proof of Lemma \ref{lem:j4}, we obtain 
$\sum a_{\mathfrak{a}} \beta_{11}(\mathfrak{a}) \ll A(x) (\log x)^{-C}$
from an application of postulate (\ref{it:biax}) with $D=D_0 D_1$,
$\ell = w$ or $\ell = x w^{-1}$ (in succession),
\[
b(\mathfrak{a}) = \mu(\max(x y^{-2}, K) \leq 
\mathfrak{a} < \min(w, K (1 + (\log N)^{-C})),\;\;\;
c(\mathfrak{a}) = 
\mathop{\mathop{\sum_{\mathfrak{d}|\mathfrak{a}}}_{
w/K < \mathfrak{d} \leq x w^{-1}/K}}_{\mathfrak{p} | \mathfrak{a}/
\mathfrak{d} \Rightarrow \mathfrak{p}\notin \mathscr{Q}} \lambda_{\mathfrak{d}}
,\]
and $K$ between $x y^{-2}$ and $w$; cf. (\ref{eq:poeque}). It remains to bound
the contribution of $\beta_{10}$, viz.,
\[\sum_{\mathfrak{a}\leq x} a_{\mathfrak{a}}
\beta_{10}(\mathfrak{a})
= 
\sum_{w<\mathfrak{b}\leq x w^{-1}}
\sum_{\mathfrak{d}|\mathfrak{b}}
\lambda_{\mathfrak{d}} \mu(\mathfrak{b})
\mathop{\sum_{\mathfrak{e}\leq x/ N \mathfrak{b}}}_{
\mathfrak{p} |\mathfrak{e} \Rightarrow \mathfrak{p} \notin
\mathscr{Q}}
\sum_{\mathfrak{c}|\mathfrak{e}}
\mu(w<\mathfrak{c}\leq x w^{-1})
a_{\mathfrak{b} \mathfrak{e}} .\]
Since (cf. (\ref{eq:rust}))
\[\begin{aligned}
\mu(w < \mathfrak{c} \leq x w^{-1}) &= \sum_{\mathfrak{c} | \mathfrak{e}}
\mu(\mathfrak{c} > w) - \sum_{\mathfrak{c} | \mathfrak{e}} 
\mu(\mathfrak{c} > x w^{-1})  \\ &= 
\mu(\rad(\mathfrak{e})) 
\sum_{\mathfrak{c} | \mathfrak{e}} \mu(\mathfrak{c} < \rad(\mathfrak{e})/w)
- \mu(\rad(\mathfrak{e}))
\sum_{\mathfrak{c} | \mathfrak{e}} \mu(\mathfrak{c} < \rad(\mathfrak{e})/x w^{-1})
,
\end{aligned}\]
we may proceed as in the proof of Lemma \ref{lem:j7} and obtain
\[\sum_{\mathfrak{a} \leq x} a_{\mathfrak{a}} \beta_{11}(\mathfrak{a})
\ll 
\frac{(\log \log x)^4 (\log \log \log x)^{\epsilon}}{
\log x} A(x) .\]
\end{proof}
\begin{prop}\label{prop:mendig}
Let $\{a_{\mathfrak{a}}\}$ be a sequence satisfying the postulates for
$n = x$. Then, for any $\epsilon>0$,
\[\sum_{\mathfrak{a}\leq x} a_{\mathfrak{a}} \mu(\mathfrak{a})
\ll_{\epsilon} \frac{(\log \log x)^4 (\log \log \log x)^{\epsilon}}{
\log x} A(x) .\]
\end{prop}
\begin{proof}
Immediate from Lemmas \ref{lem:ontont} -- \ref{lem:j4}.
\end{proof}
\section{Conclusion}\label{sec:finrem}
The Liouville function $\lambda: \mathbb{Z} - \{0\} \mapsto \{-1,1\}$
is defined as follows:
\[\lambda(n) = \prod_{p|n} (-1)^{v_p(n)} .\]
\begin{main*}
Let $f\in \mathbb{Z}\lbrack x,y\rbrack$ be an irreducible homogeneous
polynomial of degree $3$. Let $N$ be a positive number,
$S\subset \lbrack - N, N\rbrack^2$ a convex set with
$\Area(S) > N^2 (\log N)^{-A}$, and
$L$ a lattice coset of index $\lbrack \mathbb{Z}^2 : L\rbrack \ll
(\log n)^{A}$. Let $\alpha = \lambda$, $\alpha = \mu$, or
$\alpha(n) = (-1)^{\omega(n)}$. Then, for every $\epsilon>0$,
\begin{equation}\label{eq:hj}\begin{aligned}
\mathop{\sum_{(x,y)\in S\cap L}}_{\gcd(x,y)=1} \alpha(f(x,y)) 
&\ll_{f,A,\epsilon} \frac{(\log \log N)^4 (\log \log \log N)^{\epsilon}}{
\log N} \cdot \mathop{\sum_{(x,y)\in S\cap L}}_{\gcd(x,y)=1} 1 ,\\
\mathop{\sum_{(x,y)\in S\cap L}}_{(x,y)\ne (0,0)} \alpha(f(x,y)) 
&\ll_{f,A,\epsilon} \frac{(\log \log N)^4 (\log \log \log N)^{\epsilon}}{
\log N} \cdot \sum_{(x,y)\in S\cap L} 1 ,\\
\end{aligned}\end{equation}
\end{main*}
\begin{proof}
By \cite{HBM}, Lemma 2.1, there are a cubic extension $K/\mathbb{Q}$
and $\mathbb{Q}$--linearly independent elements $\omega_1, \omega_2 \in K$
such that $f(x,y) = \frac{N (x \omega_1 + y \omega_2)}{N \mathfrak{d}}$,
where $\mathfrak{d}$ is the minimal ideal of $\mathscr{O}_K$ containing
$\mathbb{Z} \omega_1 + \mathbb{Z} \omega_2$. Proposition \ref{prop:mendig}
implies that
\[\mathop{\mathop{\sum_{(x,y)\in S\cap L}}_{\gcd(x,y)=1}}_{
\mathfrak{b} | x \omega_1 + y \omega_2}
\mu\left(\frac{N (x \omega_1 + y \omega_2)}{N \mathfrak{b}}\right)
\ll_{f,\varpi,\epsilon}
\frac{(\log \log N)^4 (\log \log \log N)^{\epsilon}}{
\log N} \cdot \mathop{\mathop{\sum_{(x,y)\in S\cap L}}_{\gcd(x,y)=1}}_{
\mathfrak{b} | x \omega_1 + y\omega_2} 1 \]
for all non-zero 
ideals $\mathfrak{b}\subset \mathfrak{d}$ with 
$N \mathfrak{b} \ll (\log N)^{\varpi}$, $\varpi>A$. There
are few pairs $(x,y)$ such that $d^2|f(x,y)$ for $d$ large;
see postulate (\ref{it:axsq}), or, ultimately, \cite{Grs}. 
The rest is routine. Use, e.g., Proposition 3.2
 in \cite{Hesq} as in the proofs of Propositions 3.11--3.12
(ibid.).
Once the first equation in (\ref{eq:hj}) is proved, the second
one follows easily.
\end{proof}
\section{Final remarks}
The main theorem holds for all homogeneous cubic polynomials
$f\in \mathbb{Z}\lbrack x,y\rbrack$
irreducible in $\mathbb{Q}\lbrack x,y\rbrack$. If $\alpha = \lambda$,
this follows trivially from the theorem applied to
$\frac{1}{d_0} f$, where $d_0$ is the g.c.d. of the coefficients of $f$;
if $\alpha = \mu$ or $\alpha = (-1)^{\omega(n)}$, apply the theorem
to $\frac{1}{d_0} f$ on lattices $L'$ of the form $L' = \{(x,y) \in L :
d|f(x,y)\}$, where $d$ ranges across the integers $\ll (\log N)^{C}$
satisfying $p|d \Rightarrow p|d_0$. Thus the main theorem here and
Theorems 3.3 and 4.2 in \cite{Hered} cover together all homogeneous polynomials
$f\in\mathbb{Z}\lbrack x,y\rbrack$ of degree $3$.

Knowing that (\ref{eq:honeysuckle}) holds for $\deg f = 3$ allows us to
conclude that in certain one-parameter families of elliptic curves the
root number $W(E) = \pm 1$ averages to $0$ (\cite{Hell}, Theorem 1.2).
For example, the family $\mathscr{E}(t)$ given by 
\[c_4 = 1 - 1728 (t^3+1),\;\;\;
c_6 = (1 - 1728 (t^3+1))^2 \]
has $\av_{t\in \mathbb{Q}} W(\mathscr{E}(t)) = 0$ by the Main Theorem,
applied to the polynomial $f = x^3 + 2 y^3$. Here, as in general in
\cite{Hell}, we average over
$\mathbb{Q}$ after ordering the rationals by height.

\affiliationone{H. A. Helfgott\\
D\'epartement de Math\'ematiques et Statistique\\
Universit\'e de Montr\'eal\\CP 6128 succ Centre-Ville\\ 
Montr\'eal, QC H3C 3J7\\ Canada
\email{helfgott@dms.umontreal.ca}}
\end{document}